\documentclass{article}

\usepackage{amssymb}
\usepackage{amscd}
\usepackage{amsthm}

\newtheorem{thm}{Theorem}[section]
\newtheorem{prop}[thm]{Proposition}
\newtheorem{lem}[thm]{Lemma}
\newtheorem{cor}[thm]{Corollary}

\theoremstyle{definition}
\newtheorem{dfn}[thm]{Definition}

\theoremstyle{remark}
\newtheorem*{rem}{Remark}

\newcommand{\g}{\mathfrak{g}}
\newcommand{\Z}{\mathbb{Z}}
\newcommand{\R}{\mathbb{R}}
\newcommand{\C}{\mathbb{C}}
\newcommand{\T}{\mathbb{T}}
\newcommand{\F}{{\mathcal F}}
\newcommand{\G}{{\mathcal G}}
\newcommand{\U}{{\mathfrak U}}
\newcommand{\Ad}{\textrm{Ad}}
\newcommand{\Cech}{$\check{\textrm{C}}$ech {}}
\newcommand{\Lie}{\textrm{Lie}}
\newcommand{\Tr}{\textrm{Tr}}
\newcommand{\ev}{\textrm{ev}}
\newcommand{\id}{\textrm{id}}

\renewcommand{\i}{\sqrt{-1}}
\renewcommand{\l}{\ell}

\def\chk#1{ \check{#1} }
\def\h#1{ \widehat{#1} }
\def\til#1{ \tilde{#1} }
\def\un#1{ \underline{#1} }
\def\ol#1{ \bar{#1} }
\def\M#1#2{ C^{\infty}(#1, #2) } 

\def\hB{\widehat{B}}
\def\hb{\widehat{b}}
\def\hq{\widehat{q}}
\def\hgamma{\widehat{\gamma}}
\def\hGamma{\widehat{\Gamma}}

\title{Connections and curvings on lifting bundle gerbes}

\author{Kiyonori Gomi
\thanks{The author's research is supported by Research Fellowship of the Japan Society for the Promotion of Science for Young Scientists.}}

\date{}

\begin{document}

\maketitle

\begin{abstract}
We construct a connection and a curving on a bundle gerbe associated with lifting a structure group of a principal bundle to a central extension. The construction is based on certain structures on the bundle, i.e.\ connections and splittings. The Deligne cohomology class of the lifting bundle gerbe with the connection and with the curving coincides with the obstruction class of the lifting problem with these structures.
\end{abstract}


\section{Introduction}

Lifting the structure group of a principal bundle over $X$ to a central extension by the unit circle $\T = \{ u \in \C |\ |u| = 1 \}$ has an obstruction described by a cohomology class in $H^2(X, \un{\T}) \cong H^3(X, \Z)$. We can represent the class using the notion of bundle gerbe invented by Murray \cite{Mu}. A lifting problem defines a so-called lifting bundle gerbe whose Dixmier-Douady class is the obstruction class \cite{Mu, Mu-S1}.

In \cite{B}, Brylinski defined a class in the Deligne cohomology group
$$
H^2(X, \un{\T} \to \i \un{\Omega}^1 \to \i \un{\Omega}^2)
$$
which expresses an obstruction of the lifting problem with certain geometric structures: a connection and a \textit{splitting}. A splitting is formulated as a map of vector bundles associated with the underlying bundle. Using a splitting, we define the \textit{scalar curvature} as the `scalar component' of the curvature of a connection on a lifting compatible with the connection given. The Deligne cohomology class is the obstruction to the construction of a lifting with a compatible connection whose scalar curvature vanishes.

Any bundle gerbe admits differential geometric structures called \textit{connections} and \textit{curvings (bundle gerbe curvatures)} \cite{Mu, Mu-S1}. A bundle gerbe with such structures defines a class in the Deligne cohomology group. This correspondence gives a bijection from the stable isomorphism classes of bundle gerbes with these structures to the Deligne cohomology group.

In this paper we construct explicitly a connection and a curving on a lifting bundle gerbe which give the same Deligne cohomology class as the obstruction class. For the purpose we choose a split of the Lie algebra of the central extension as a vector space, and use a notion of \textit{reduced splitting}. If a split of the Lie algebra is fixed, then splittings correspond to reduced splittings bijectively.

\smallskip

Note that a class in $H^3(X, \Z)$ is also represented by the other geometric object called \textit{gerbe} \cite{B,B-M,Gi}. As an original example of gerbes, the lifting problem was investigated by Brylinski in \cite{B}. In a natural fashion  a connection and a splitting induce differential geometric structures on the gerbe associated with the problem. The gerbe with these structures also provides the obstruction class of lifting in the Deligne cohomology group.

\smallskip

As an example we study the loop bundle \cite{C-P, K} $LP \to LM$ induced from a principal $SU(2)$-bundle $P \to M$. In this case a connection $A$ on $P$ defines both a connection and a reduced splitting of $LP$. The curving coincides with the 2-form on $LP$ defined by Coquereaux and Pilch \cite{C-P} in extending a relation between 2-forms on the loop group to a relation on $LP$. As a consequence we can relate the curving to the transgression of the Chern-Simons form \cite{C-S}, and the 3-curvature \cite{Mu, Mu-S1} to the transgression of the characteristic form defined by the connection on $P$. In addition we study the case that we have a compact oriented 2-manifold $\Sigma$ with boundary $\partial \Sigma = S^1$. We construct a lifting of the pull-back of $LP$ under the restriction map $r : \M{\Sigma}{M} \to LM$. We give a compatible connection, and compute the scalar curvature using the curving.

\smallskip

Recently, Murray and Stevenson \cite{Mu-S2} also described a connection and a curving on the lifting bundle gerbe associated with a loop group bundle. They use \textit{twisted Higgs fields} to define curvings. We see that both expressions of the reduced splitting and of the curving discussed in our example include twisted Higgs fields.

\medskip

The organization of this paper is as follows. In Section \ref{sec_deligne} we introduce Deligne cohomology groups. In Section \ref{sec_lift_prob} we first recall ordinary lifting problem. We next introduce splitting and scalar curvature, and explain the problem with additional structures. We also introduce the notion of reduced splitting here. In Section \ref{sec_bundle_gerbe} we overview bundle gerbes. Section \ref{sec_construction} contains the main result. We construct the connection and the curving on a lifting bundle gerbe. In Section \ref{sec_example} we study an example of lifting problems associated with loop bundles.


\section{Deligne cohomology groups}
\label{sec_deligne}

We introduce Deligne cohomology groups which are certain refinements of usual cohomology groups.

\begin{dfn}[\cite{B}]
Let $n, m$ be non-negative integers and $X$ a smooth manifold. We denote the sheaf of $\T$-valued functions by $\un{\T}$ and the sheaf of $\R$-valued differential $q$-forms by $\un{\Omega}^q$. We denote by $\F^m$ the following complex of sheaves
$$
\begin{CD}
\un{\T} @>{d \log}>> \i \un{\Omega}^1 @>d>> \ldots @>d>> \i \un{\Omega}^m.
\end{CD}
$$
We define the \textit{Deligne cohomology group} by the hypercohomology group of the complex of sheaves $H^n(X, \F^m)$.
\end{dfn}

Usually hypercohomology groups are computed by \Cech cohomology groups \cite{B-T, B}. If we take an open cover $\U$ of $X$, then we have a double complex. One coboundary operator comes from that on $\F^m$. Another operator is the \Cech coboundary operator. We denote by $\chk{H}^n(\U, \F^m)$ the cohomology of the total complex. The hypercohomology $H^n(X, \F^m)$ is isomorphic to the direct limit of the \Cech cohomology groups with respect to refinements of open covers. If we take $\U$ as a $good$ cover \cite{B-T, B}, i.e.\ all the intersections of a finite number of open sets belonging to $\U$ are contractible, then we can naturally identify the \Cech cohomology group $\chk{H}^n(\U, \F^m)$ with the hypercohomology group $H^n(X, \F^m)$.

\begin{thm}[\cite{B}]
We have an exact sequence
$$
\begin{CD}
0 @>>> H^n(X, \T) @>>> H^n(X, \F^n) @>d>> 2\pi\i \Omega^{n+1}(X)_0 @>>>0,
\end{CD}
$$
where $\Omega^{n+1}(X)_0$ is the group of closed $(n+1)$-forms on $X$ whose periods are integers. The map $d$ is defined by $d(g, \omega^1, \ldots, \omega^n) =  d\omega^n$ on the \Cech cocycle.
\end{thm}


\section{Lifting problems}
\label{sec_lift_prob}

In this section we explain the lifting problems and the associated obstruction classes. We mainly refer the book of Brylinski \cite{B}. However, the formulations here are slightly different from those in the book.

\smallskip

Let $\Gamma$ be a Lie group. We fix a central extension $\hGamma$ of $\Gamma$ by the unit circle $\T = \{ u \in \C |\ |u| = 1 \}$
\begin{eqnarray}
\begin{CD}
1 @>>> \T @>>> \hGamma @>q>> \Gamma @>>> 1.
\end{CD}
\label{es_ce_group}
\end{eqnarray}
This induces a central extension of the Lie algebra $\Lie\Gamma$ by $\Lie\T = \i\R$
\begin{eqnarray}
\begin{CD}
0 @>>> \i \R @>>> \Lie \hGamma @>q_*>> \Lie \Gamma @>>> 0.
\end{CD}
\label{es_ce_liealg}
\end{eqnarray}
We always identify the elements contained in the center of $\hGamma$ and $\Lie\hGamma$ with $\T$ and $\i\R$ respectively without mentioning. The Lie brackets of these Lie algebras are denoted as $[\h{X}, \h{Y}]_{\hGamma}$ and $[X, Y]_{\Gamma}$.

If we replace the center $\T$ by $\C^*$, then we obtain corresponding results under appropriate modifications.

\begin{dfn}[\cite{B}]
Let $\pi : B \to X$ be a principal $\Gamma$-bundle. A \textit{lifting} $(\hB, \hq)$ is a $\hGamma$-bundle $\h{\pi} : \hB \to X$ with an equivariant map $\hq : \hB \to B$ which satisfies $\hq(\hb \cdot \hgamma) = \hq(\hb) \cdot q(\hgamma)$ for $\hb \in \hB$ and $\hgamma \in \hGamma$. 
\end{dfn}

\begin{thm}[\cite{B}] \label{thm_lp_top}
Let $\{ g_{\alpha \beta} \}$ be the transition functions of a $\Gamma$-bundle $B \to X$ defined by local sections $\{ s_\alpha \}$ with respect to a good cover $\U$. Take $\h{g}_{\alpha \beta} : U_\alpha \cap U_\beta \to \hGamma$ such that $q(\h{g}_{\alpha \beta}) = g_{\alpha \beta}$, and define $z_{\alpha \beta \gamma} : U_\alpha \cap U_\beta \cap U_\gamma \to \T$ by 
\begin{eqnarray}
z_{\alpha \beta \gamma} =
\h{g}_{\alpha \beta} \h{g}_{\beta \gamma} \h{g}_{\gamma \alpha}.
\label{formula_dfn_z}
\end{eqnarray}

{\bf (a)} The \Cech cochain $(z_{\alpha \beta \gamma}) \in \chk{C}^2(\U, \un{\T})$ defines a cohomology class in $H^2(X, \un{\T})$ which depends only on $B$.

{\bf (b)} There is a lifting of $B$ if and only if the class is trivial in $H^2(X, \un{\T})$.
\end{thm}

The proof can be found in the proof of Theorem \ref{thm_lp_con_spl} and is omitted.

\medskip

Next we consider lifting problems with additional structures.

\begin{dfn}[\cite{B}]
Let $\theta \in \Omega^1(B; \Lie \Gamma)$ be a connection on a $\Gamma$-bundle $B \to X$. Suppose that there exists a lifting $(\hB, \hq)$ of $B$. A connection $\h{\theta}$ compatible with $\theta$ is a connection on $\hB$ such that $q_*(\h{\theta}) = \hq^* \theta$.
\end{dfn}

\begin{lem}[\cite{B}] \label{lem_cpb_conn}
The set of connections compatible with $\theta$ is an affine space under the vector space $\Omega^1(X; \i \R)$.
\end{lem}

\begin{proof}
If $\h{\theta}$ and $\h{\theta}'$ are connections compatible with $\theta$, then the difference $\h{\theta}' - \h{\theta}$ is a 1-form on $\hB$ with values in $\i\R$, the center of the Lie algebra of $\hGamma$. By the property of connection, the 1-form on $\hB$ descends to a 1-form on $X$. Conversely, if $\alpha$ is a 1-form on $X$ with values in $\i\R$, then $\h{\theta} + \h{\pi}^*\alpha$ is a connection on $\hB$. Clearly this connection is compatible with $\theta$.
\end{proof}

The Lie group $\Gamma$ acts on $\Lie \Gamma$ by the adjoint action. We can lift this action to that on $\Lie \hGamma$ as follows. For $\gamma \in \Gamma$ we take an element $\hgamma \in \hGamma$ such that $q(\hgamma) = \gamma$. Using the adjoint action of $\hGamma$ on $\Lie \hGamma$, we define $\Ad_\gamma : \Lie\hGamma \to \Lie\hGamma$ by $\Ad_{\gamma}\h{X} = \Ad_{\hgamma}\h{X}$. Since the center $\T$ acts on $\Lie\hGamma$ trivially, we have $\Ad_{\hgamma}\h{X} = \Ad_{u\hgamma} \h{X}$ for all $u \in \T$ and $\Ad_\gamma \h{X}$ is well-defined.

As a result the Lie group $\Gamma$ acts on the exact sequence (\ref{es_ce_liealg}), and we obtain an exact sequence of the vector bundles associated with a $\Gamma$-bundle $B$
\begin{eqnarray}
\begin{CD}
0 @>>> B \times_{Ad} (\i \R) \\
@>>> B \times_{Ad} \Lie \hGamma
@>>> B \times_{Ad} \Lie \Gamma @>>>0.
\end{CD}
\end{eqnarray}
This implies that $B \times_{Ad} \Lie \hGamma$ is isomorphic to the direct sum of $B \times_{Ad} \Lie\Gamma$ and $B \times_{Ad} (\i \R) \cong X \times (\i \R)$ as vector bundles. A bundle map which gives this isomorphism is a \textit{splitting}.

\begin{dfn}[\cite{B}]
A \textit{splitting} of $B$ is a vector bundle map
\begin{eqnarray}
L : B \times_{Ad} \Lie \hGamma \rightarrow B \times_{Ad} (\i \R)
\end{eqnarray}
which is identity on the subbundle $B \times_{Ad} (\i \R)$.
\end{dfn}

\begin{prop}
The set of splittings of $B$ is an affine space under the vector space of sections of the vector bundle $\textrm{Hom}(B \times_{Ad} \Lie \Gamma, X \times {\i \R})$ over $X$.
\end{prop}

\begin{proof}
For a splitting $L$ we put $L([b, \h{X}]) = [b, \ol{L}(b, \h{X})]$. If splittings $L$ and $L'$ are given, then we have a map $\varepsilon : B \times \Lie\hGamma \to \i\R$ defined by $\varepsilon(b, \h{X}) = \ol{L}'(b, \h{X}) - \ol{L}(b, \h{X})$. It is easy to see that $\varepsilon$ is trivial on the subspace $B \times \i\R$ and satisfies $\varepsilon(b, \h{X}) = \varepsilon(b\gamma, \Ad_{\gamma^{-1}}\h{X})$. Hence $\varepsilon$ induces a map $B \times \Lie\Gamma \to \i\R$ and gives rise to a section of $\textrm{Hom}(B \times_{Ad} \Lie\Gamma, B \times_{Ad} (\i\R))$. Conversely, a section of the vector bundle produces a new splitting by the translation.
\end{proof}

By the help of a splitting, we can extract the ``scalar component" from the curvature of a compatible connection. Note that for any lifting $(\hB, \hq)$ of $B$ we have a natural isomorphism $\hB \times_{Ad} \Lie\hGamma \to B \times_{Ad} \Lie\hGamma$ of vector bundles defined by $[\hb, \h{X}] \mapsto [\hq(\hb), \h{X}]$. 

\begin{dfn}[\cite{B}]
Let $\hB$ be a lifting of a $\Gamma$-bundle $B$ over $X$, $\h{\theta}$ a connection on $\hB$, and $L$ a splitting of $B$. We regard the curvature $F_{\h{\theta}} = d\h{\theta} + \frac{1}{2}[\h{\theta}, \h{\theta}]_{\hGamma}$ as a 2-form on $X$ which takes its values in $\hB \times_{Ad} \Lie\hGamma \cong B \times_{Ad} \Lie\hGamma$. We define the \textit{scalar curvature} of $\h{\theta}$ as a $\i\R$-valued 2-form on $X$ obtained by composing the splitting:
\begin{eqnarray}
K_{\h{\theta}} = L \circ F_{\h{\theta}}.
\end{eqnarray}
\end{dfn}

\begin{lem}[\cite{B}]
Let $\h{\theta}$ be a connection on a lifting $\hB$, and $\h{\pi} : \hB \to X$ the projection. For $\alpha \in \i\Omega^1(X)$ the scalar curvature of $\h{\theta} + \h{\pi}^*\alpha$ is
\begin{eqnarray}
K_{\h{\theta} + \h{\pi}^*\alpha} = K_{\h{\theta}} + d\alpha .
\end{eqnarray}
\end{lem}

\begin{proof}
Since $F_{\h{\theta} + \h{\pi}^*\alpha} = F_{\h{\theta}} + \h{\pi}^*d\alpha$ and a splitting is the identity on the subbundle $X \times \i\R$, we have the result.
\end{proof}

\begin{thm}[\cite{B}] \label{thm_lp_con_spl}
Let $B \to X$ be a $\Gamma$-bundle with a connection $\theta$ and a splitting $L$. For local sections $\{ s_\alpha \}$ with respect to a good cover $\U$, we denote the transition functions by $g_{\alpha \beta}$ and the connection forms by $\theta_\alpha = s_\alpha^*\theta$. Take $\h{g}_{\alpha \beta} : U_\alpha \cap U_\beta \to \hGamma$ such that $q(\h{g}_{\alpha \beta}) = g_{\alpha \beta}$ and $\h{\theta}_\alpha \in \Omega^1(U_\alpha; \Lie \hGamma)$ such that $q_*(\h{\theta}_\alpha) = \theta_\alpha$. We define $z_{\alpha \beta \gamma}$, $u_{\alpha \beta}$ and $K_{\alpha}$ by
\begin{eqnarray}
\left\{
\begin{array}{ccl}
z_{\alpha \beta \gamma} 
& = & 
\h{g}_{\alpha \beta} \h{g}_{\beta \gamma} \h{g}_{\gamma \alpha}, \\
u_{\alpha \beta} 
& = &
\h{\theta}_\beta -
\left\{ 
{\h{g}_{\alpha \beta}}^{-1} \h{\theta}_\alpha \h{g}_{\alpha \beta} 
+ \h{g}_{\alpha \beta}^* \h{\mu}
\right\}, \\
K_{\alpha}
& = &
\ol{L}( s_\alpha, F_{ \h{\theta}_\alpha } ),
\end{array}
\right.
\end{eqnarray}
where $\h{\mu}$ is the Maurer-Cartan form of $\hGamma$ and $\ol{L} : B \times \Lie\hGamma \to \i\R$ is defined by $L([b, \h{X}]) = [b, \ol{L}(b, \h{X})]$. 

{\bf (a)} The \Cech cochain $(z_{\alpha \beta \gamma}, u_{\alpha \beta}, K_{\alpha}) \in \chk{C}^2(\U, \F^2)$ defines a cohomology class in $H^2(X, \F^2)$ which depends only on $(B, \theta, L)$.

{\bf (b)} There exists a lifting of $B$ and a connection compatible with $\theta$ whose scalar curvature vanishes if and only if the class is trivial in $H^2(X, \F^2)$.
\end{thm}

\begin{proof}
First note that the values of $z_{\alpha \beta \gamma}$ and $u_{\alpha \beta}$ are contained in the center of $\hGamma$ and of $\Lie\hGamma$ respectively as a consequence of the cocycle condition
\begin{eqnarray}
g_{\alpha \gamma} & = & g_{\alpha \beta} g_{\beta \gamma}, \\
\theta_\beta & =&
\Ad_{ g_{\alpha \beta}^{-1} } \theta_\alpha + g_{\alpha \beta}^*\mu,
\end{eqnarray}
where $\mu$ is the Maurer-Cartan form of $\Gamma$. For (a), the following equalities obtained by direct computation show the cocycle condition of $(z_{\alpha \beta \gamma}, u_{\alpha \beta}, K_\alpha)$:
\begin{eqnarray}
d \log z_{\alpha \beta \gamma} 
& = &
\Ad_{\h{g}_{\gamma \alpha}^{-1}} \Ad_{\h{g}_{\beta \gamma}^{-1}}
\h{g}_{\alpha \beta}^*\h{\mu}
+
\Ad_{\h{g}_{\gamma \alpha}^{-1}} \h{g}_{\beta \gamma}^*\h{\mu}
+
\h{g}_{\gamma \alpha}^*\h{\mu}, \\
F_{\h{\theta}_\beta} & = & {\h{g}_{\alpha \beta}}^{-1} F_{\h{\theta}_\alpha} \h{g}_{\alpha \beta} + du_{\alpha \beta}.
\end{eqnarray}
We can verify that resulting \Cech cohomology class is independent of the choice of local section $s_\alpha$ and of the choice of $\h{g}_{\alpha \beta}$, $\h{\theta}_\alpha$. The standard argument shows that the definition of the class is independent of the choice of $\U$. For (b), if we have a lifting and a compatible connection whose scalar curvature vanishes, then we have $z_{\alpha \beta \gamma} = 1$ and $u_{\alpha \beta} = 0$ using the transition functions and the connection forms of the lifting. Because $K_\alpha$ is the local description of the scalar curvature, we also have $K_\alpha = 0$. Hence the \Cech cocycle is trivial. Conversely, if the class is the coboundary of a \Cech cochain $(h_{\alpha \beta}, k_\alpha)$, then we put
\begin{eqnarray}
\left\{
\begin{array}{ccl}
\h{g}'_{\alpha \beta} & = & {h_{\alpha \beta}}^{-1} \h{g}_{\alpha \beta}, \\
\theta'_\alpha & = & \theta_\alpha - k_\alpha.
\end{array}
\right.
\end{eqnarray}
These data give a lifting $\hB$ and a compatible connection $\h{\theta}$. The scalar curvature is calculated as $K_{\h{\theta}}|_{U_\alpha} = K_\alpha - dk_\alpha = 0$ on each $U_\alpha$. 
\end{proof}

In the remainder of this section we introduce and study an equivalent description of a splitting which we will use in Section \ref{sec_construction}.

\begin{dfn} \label{dfn_group_cocycle}
Let $\sigma : \Lie\Gamma \to \Lie\hGamma$ be a split of the exact sequence (\ref{es_ce_liealg}) as vector spaces. We define a map $Z_\sigma : \Gamma \times \Lie \Gamma \to \i \R$ by 
\begin{eqnarray}
Z_\sigma(\gamma, X) = \Ad_{\gamma}\sigma(X) - \sigma(\Ad_\gamma X).
\end{eqnarray}
\end{dfn}

We call $Z_\sigma$ the \textit{group cocycle} for the central extension $\hGamma$, because it indeed gives a group 1-cocycle \cite{Bro} of $\Gamma$ with coefficients in $\textrm{Hom}(\Lie\Gamma, \i\R)$.

\begin{lem}
For any $\gamma, \eta \in \Gamma$ and $X \in \Lie\Gamma$ we have
\begin{eqnarray}
Z_\sigma(\gamma \eta, X) = Z_\sigma(\gamma, \Ad_\eta X) + Z_\sigma(\eta, X).
\label{formula_Z_Ad}
\end{eqnarray}
If $\sigma'$ is the other split, then we have
\begin{eqnarray}
Z_{\sigma'}(\gamma, X) - Z_\sigma(\gamma, X) = 
(\sigma' - \sigma)(X - \Ad_\gamma X),
\label{formula_Z_diff_sigma}
\end{eqnarray}
where $\sigma' - \sigma$ is regarded as a linear map $\Lie\Gamma \to \i\R$.
\end{lem}

\begin{proof}
We can directly prove these formulas by using the fact that the adjoint action of $\hGamma$ on the center $\i\R$ is trivial.
\end{proof}

\begin{rem}
If a split $\sigma$ is given, then the Lie algebra cocycle for the central extension $\Lie \hGamma$ is defined by $\omega_\sigma(X, Y) = [\sigma(X), \sigma(Y)]_{\hGamma} - \sigma([X, Y]_{\Gamma})$. There is a relation between $Z_\sigma$ and $\omega_\sigma$:
$$
\frac{d}{dt} \bigg|_{t=0} Z_\sigma( \exp{tX}, Y) = \omega_\sigma(X, Y).
$$
\end{rem}

\begin{dfn}
A \textit{reduced splitting} of a $\Gamma$-bundle $B$ with respect to a split $\sigma$ is defined as a map $\l_\sigma : B \times \Lie \Gamma \to \i\R$ which is linear on each element of $B$ and satisfies
\begin{eqnarray}
\l_\sigma(b, X) = 
\l_\sigma(b \gamma, \Ad_{\gamma^{-1}}X) + Z_\sigma(\gamma^{-1}, X)
\label{formula_red_split}
\end{eqnarray}
for all $(b, X) \in B \times \Lie\Gamma$ and $\gamma \in \Gamma$.
\end{dfn}

Notice that the relation (\ref{formula_red_split}) is equivalent to
\begin{eqnarray}
\l_\sigma(b\gamma, X) =
\l_\sigma(b, \Ad_{\gamma}X) + Z_\sigma(\gamma, X).
\end{eqnarray}
Thus the cocycle $Z_\sigma$ defines a reduced splitting of the trivial bundle $X \times \Gamma$.

\begin{thm} \label{thm_spl_redspl}
Let $B$ be a $\Gamma$-bundle and $\sigma$ a split of the Lie algebra. There is a bijective correspondence between splittings of $B$ and reduced splittings of $B$ with respect to $\sigma$.
\end{thm}

\begin{proof}
We write $L([b, \h{X}]) = [b, \ol{L}(\hb, \h{X})]$ for a splitting $L$. For $X \in \Lie\Gamma$ we have
\begin{eqnarray*}
\ol{L}(b, \sigma(X)) 
& = & 
\ol{L}(b \gamma, \Ad_{\gamma^{-1}}\sigma(X))
  = 
\ol{L}(b \gamma, 
\sigma(\Ad_{\gamma^{-1}} X) + Z_\sigma(\gamma^{-1}, X)) \\
& = &
\ol{L}(b \gamma, \sigma(\Ad_{\gamma^{-1}} X))
+ Z_\sigma(\gamma^{-1}, X).
\end{eqnarray*}
Hence we get a reduced splitting with respect to $\sigma$. Conversely, let $\l_\sigma$ be a reduced splitting. We can uniquely decompose any element $\h{X} \in \Lie\h\Gamma$ as $\h{X} = \sigma(X) + z$, where $X = q_*(\h{X})$ and $z$ is contained in the center. If we put $\til{\l}_\sigma(b, \h{X}) = \l_\sigma(b, X) + z$, then $[b, \h{X}] \mapsto [b, \til{\l}_\sigma(b, \h{X})]$ is a splitting of $B$.
\end{proof}

\begin{cor} \label{cor_spl_redspl}
The splittings of $B$ induced from reduced splittings $\l_\sigma$ and $\l_{\sigma'}$ are identical if and only if we have
\begin{eqnarray}
\l_{\sigma'}(b, X) - \l_{\sigma}(b, X) = (\sigma' - \sigma)(X)
\label{rel_redspl}
\end{eqnarray}
for all $(b, X) \in B \times \Lie\Gamma$.
\end{cor}


\section{Bundle gerbes}
\label{sec_bundle_gerbe}

To make this paper self-contained we briefly explain \textit{bundle gerbes} which were originally invented by Murray. Details are in the papers \cite{Mu, Mu-S1}, and we omit the proofs. As a bundle gerbe is a geometric object representing a degree 3 integral cohomology class, a \textit{gerbe} \cite{B, B-M, Gi} is also such a geometric object. The relation between gerbes and bundle gerbes can be seen in \cite{Mu, Mu-S1} too.

\smallskip

A map $\pi : Y \to X$ is called \textit{locally split} \cite{Mu-S1} if the following condition holds: for every $x \in X$ there is an open set $U$ containing $x$ and a local section $s : U \to Y$. By the definition a locally split map is surjective. Locally trivial fiber bundles give examples of locally split maps. We define the fiber product by $Y^{[2]} = \{(y_1, y_2) \in Y \times Y | \pi(y_1) = \pi(y_2) \}$. The $p$-fold fiber product $Y^{[p]}$ is similarly defined. We define the projection $\pi_{1 2 \ldots \hat{i} \ldots p} : Y^{[p]} \to Y^{[p-1]}$ by omitting $i$-th component.

\begin{dfn}[\cite{Mu, Mu-S1}]
Let $(P, Y)$ be a pair consisting of a locally split map $\pi : Y \to X$ and a principal $\T$-bundle $P$ over $Y^{[2]}$. A \textit{product} on $P$ is defined by a bundle isomorphism $m : \pi_{12}^*P \otimes \pi_{23}^*P \to \pi_{13}^*P$ which is associative whenever the triple product exists. A pair $(P, Y)$ together with a product on $P$ is called a \textit{bundle gerbe} over $X$.
\end{dfn}

An isomorphism of bundle gerbes $(P, Y)$ and $(P', Y')$ is a pair $(\til{F}, F)$ where $F : Y \to Y'$ is a fiber preserving map and $\til{F} : P \to P'$ is an isomorphism of $\T$-bundles which covers the induced map $F^{[2]}$ and commutes with the products on $P$ and $P'$. 

A $\T$-bundle $Q$ over $Y$ defines a bundle gerbe $(\delta(Q), Y)$ by $\delta(Q) = \pi_1^*Q^* \otimes \pi_2^*Q$, where $Q^*$ is the inverse of the $\T$-bundle $Q$. A bundle gerbe is called \textit{trivial} if it is isomorphic to $(\delta(Q), Y)$ for some $Q$ and $Y$.

There is a characteristic class for bundle gerbes defined as follows. Let $(P, Y)$ be a bundle gerbe over $X$, and $\U$ a sufficiently fine good cover. We abbreviate the notation of intersections as $U _{\alpha_1 \ldots \alpha_N} = U_{\alpha_1} \cap \cdots \cap U_{\alpha_N}$. Since $\pi : Y \to X$ is locally split, we can take local sections $s_\alpha : U_\alpha \to X$. These sections induce $(s_\alpha, s_\beta) : U_{\alpha \beta} \to Y^{[2]}$. Because $U_{\alpha \beta}$ is contractible there is a section $\sigma_{\alpha \beta} : U_{\alpha \beta} \to (s_\alpha, s_\beta)^*P$. Using the product on $P$ we have $g_{\alpha \beta \gamma} : U_{\alpha \beta \gamma} \to \T$ defined by $m(\sigma_{\alpha \beta}, \sigma_{\beta \gamma}) = \sigma_{\alpha \gamma} g_{\alpha \beta \gamma}$. Then $(g_{\alpha \beta \gamma})$ defines a cohomology class $d(P) = d(P, Y) \in H^2(X, \un{\T}) \cong H^3(X, \Z)$ called the \textit{Dixmier-Douady class}.

If a bundle gerbe $(P, Y)$ over X and a map $\varphi : X' \to X$ are given, then we have the pull-back of the bundle gerbe $\varphi^*(P, Y) = ((\til{\varphi}^{[2]})^*P, \varphi^*Y)$ over $X'$, where $\til{\varphi} : \varphi^*Y \to Y$ is the map covering $\varphi$. The Dixmier-Douady class is natural in the sense that $\varphi^*d(P, Y) = d(\varphi^*(P, Y))$.

\begin{thm}[\cite{Mu, Mu-S1}]
A bundle gerbe $(P, Y)$ over $X$ is trivial if and only if the Dixmier-Douady class is trivial in $H^2(X, \un{\T}) \cong H^3(X, \Z)$.
\end{thm}

It is known that we can construct bundle gerbes which are not isomorphic but have the same Dixmier-Douady class. Hence the notion of \textit{stable isomorphism} is introduced in \cite{Mu-S1}. The product of bundle gerbes $(P, Y)$ and $(P', Y')$ is defined by $(P \otimes P', Y \times_\pi Y')$. Bundle gerbes $(P, Y)$ and $(Q, Z)$ are stably isomorphic if there exist trivial bundle gerbes $T_1$ and $T_2$ such that the bundle gerbes $P \otimes T_1$ and $Q \otimes T_2$ are isomorphic.

\begin{thm}[\cite{Mu-S1}] \label{thm_classify_bg}
There is a bijective correspondence between the stable isomorphism classes of bundle gerbes over $X$ and $H^2(X, \un{\T}) \cong H^3(X, \Z)$.
\end{thm}

Next we introduce differential geometric structures on bundle gerbes.

\begin{dfn}[\cite{Mu, Mu-S1}]
Let $(P, Y)$ be a bundle gerbe over $X$.

{\bf (a)}
A connection $\nabla$ on the $\T$-bundle $P$ is called a \textit{bundle gerbe connection}, or a \textit{connection}, if it is compatible with the product: $\pi_{12}^*\nabla + \pi_{23}^*\nabla = m^*\pi_{13}^*\nabla$.

{\bf (b)}
A \textit{curving} for a bundle gerbe connection $\nabla$ is a 2-form $f \in \i\Omega^2(Y)$ such that $F_\nabla = - \pi_1^*f + \pi_2^*f$, where $F_\nabla$ is the curvature of the connection $\nabla$ on the $\T$-bundle $P$.
\end{dfn}

There are a bundle gerbe connection for any bundle gerbe, and a curving for any bundle gerbe connection. Note that the choice of bundle gerbe connections and curvings is not unique.

Let $Q \to Y$ be a $\T$-bundle. A connection $\alpha$ on $Q$ defines a bundle gerbe connection $\delta(\alpha) = -\pi_1^*\alpha + \pi_2^*\alpha$ on $\delta(Q)$. For the bundle gerbe connection we can take $F_\alpha$ as a curving. A bundle gerbe with connection and curving is called trivial if it is isomorphic to $(\delta(Q), \delta(\alpha), F_\alpha)$ for some $Q$ and $\alpha$.

For a bundle gerbe with connection and curving $(P, \nabla, f)$ we define a class in the Deligne cohomology group $H^2(X, \F^2)$. On $(s_\alpha, s_\beta)^*P$ we have the induced connection $\nabla_{\alpha \beta} = (s_\alpha, s_\beta)^*\nabla$. Using the section $\sigma_{\alpha \beta}$, we put $A_{\alpha \beta} = - \sigma_{\alpha \beta}^*\nabla_{\alpha \beta}$. We also put $f_\alpha = - s_\alpha^*f$. Then $(g_{\alpha \beta \gamma}, A_{\alpha \beta}, f_\alpha)$ defines a class in $H^2(X, \F^2)$. The pull-back of bundle gerbe with connection and curving is defined by $\varphi^*(P, \nabla, f) = ( (\til{\varphi}^{[2]})^*P, (\til{\varphi}^{[2]})^*\nabla, \til{\varphi}^*f)$. The Deligne cohomology class also satisfies the naturality under the pull-back.

\begin{thm}[\cite{Mu-S1}] \label{thm_triv_del_bgcc}
A bundle gerbe with connection and curving is trivial if and only if the corresponding class $[(g_{\alpha \beta \gamma}, A_{\alpha \beta}, f_\alpha)]$ is trivial in $H^2(X, \F^2)$.
\end{thm}

The notion of stable isomorphism is defined similarly.

\begin{thm}[\cite{Mu-S1}] \label{thm_classify_bgcc}
There is a bijective correspondence from the stable isomorphism classes of bundle gerbes with connection and curving over $X$ to the Deligne cohomology group $H^2(X, \F^2)$.
\end{thm}

For a curving $f$ of a connection $\nabla$ on $(P, Y)$ there exists $\Xi \in \i\Omega^3(X)$ such that $df = \pi^*\Xi$. This 3-form is called the \textit{3-curvature} of $(P, \nabla, f)$. It is known that $\frac{1}{2\pi\i}\Xi$ is a de Rham representative of $d(P) \in H^3(X, \Z)$. The 3-curvature is computed from the Deligne cohomology class $[(g_{\alpha \beta \gamma}, A_{\alpha \beta}, f_\alpha)]$ as $\Xi|_{U_\alpha} = -df_\alpha$.


\section{Main results}
\label{sec_construction}

First, we recall the \textit{lifting bundle gerbe} associated with an ordinary lifting problem. We fix a central extension $\hGamma$ of a Lie group $\Gamma$.

\begin{dfn}[\cite{Mu, Mu-S1}]
Let $\pi : B \to X$ be a $\Gamma$-bundle. We define a map $\zeta : B^{[2]} \to \Gamma$ by $b_1 \zeta(b_1, b_2) = b_2$. On the pull-back $J_B = \zeta^*\hGamma$ there is a product $m : \pi_{12}^*J_B \otimes \pi_{23}^*J_B \to \pi_{13}^*J_B$ induced by the group product in $\hGamma$. This bundle gerbe $(J_B, B)$ over $X$ is called the \textit{lifting bundle gerbe}.
\end{dfn}

\begin{thm}[\cite{Mu}] \label{thm_corres1}
The Dixmier-Douady class $d(J_B, B) \in H^2(X, \un{\T})$ of the lifting bundle gerbe associated with a lifting problem of a $\Gamma$-bundle $B$ coincides with the obstruction class of the lifting problem.
\end{thm}

\begin{proof}
Note that $(s_\alpha, s_\beta)^*J_B = g_{\alpha \beta}^*\hGamma$, where $g_{\alpha \beta}$ is the transition function of $B$. If we take the section as $\sigma_{\alpha \beta} = \h{g}_{\alpha \beta}$, then we have $g_{\alpha \beta \gamma} = z_{\alpha \beta \gamma}$.
\end{proof}

\begin{cor}[\cite{Mu}] \label{cor_top_lbg}
A $\Gamma$-bundle $B$ admits a lifting if and only if the lifting bundle gerbe $(J_B, B)$ is trivial.
\end{cor}

Now we construct a connection on the lifting bundle gerbe. For the purpose we take a split $\sigma$ of the Lie algebra, and use a connection on $B$ and the group cocycle $Z_\sigma$ of $\Gamma$.

\begin{lem}
Let $\h{\mu}$ and $\mu$ be the Maurer-Cartan forms on $\hGamma$ and $\Gamma$ respectively. If a split $\sigma : \Lie\Gamma \to \Lie\hGamma$ is given, then a connection on the $\T$-bundle $q : \hGamma \to \Gamma$ is obtained by $\nu_\sigma = \h{\mu} - \sigma(q^*\mu)$. The curvature of the connection is $F_{\nu_\sigma} = -\frac{1}{2} \omega_\sigma(\mu, \mu)$, where $\omega_\sigma$ is the Lie algebra cocycle for the central extension $\Lie\hGamma$ defined by $\omega_\sigma(X, Y) = [\sigma(X), \sigma(Y)]_{\hGamma} - \sigma([X, Y]_{\Gamma})$.
\end{lem}

\begin{proof}
Since we have $q_*(\h{\mu}) = \hq^*\mu$, the 1-form $\nu_\sigma$ takes values in $\i\R$. The basic properties of the Maurer-Cartan form show that the 1-form is a connection on the $\T$-bundle. We can easily compute the curvature of this connection using the Maurer-Cartan equation.
\end{proof}

\begin{lem}
{\bf (a)}
For $\hgamma \in \hGamma$ we have
\begin{eqnarray}
L_{\hgamma}^*\nu_\sigma & = & \nu_\sigma, \label{formula_left_nu}\\
R_{\hgamma}^*\nu_\sigma & = & \nu_\sigma + Z_\sigma(q(\hgamma)^{-1}, q^*\mu). 
\label{formula_right_nu}
\end{eqnarray}

{\bf (b)}
For a connection $\theta$ on $B$ we have
\begin{eqnarray}
\zeta^*\mu = - \Ad_{\zeta^{-1}}\pi_1^*\theta + \pi_2^*\theta .
\label{formula_mu_theta}
\end{eqnarray}
\end{lem}

\begin{proof}
The formulas (\ref{formula_left_nu}) and (\ref{formula_right_nu}) are directly proved. The formula (\ref{formula_mu_theta}) is proved as follows. A tangent vector at $(b_1, b_2) \in B^{[2]}$ is a pair $(V_1, V_2) \in T_{b_1}B \oplus T_{b_2}B$ such that $\pi_*V_1 = \pi_*V_2$. The push-forward of the vector under $\zeta$ can be expressed as $\zeta_*(V_1, V_2) = {L_{\zeta(b_1, b_2)}}_* X_{12}$, where $X_{12} \in \Lie\Gamma = T_e\Gamma$. Then we have an equality of tangent vectors $V_2 - {R_{\zeta(b_1, b_2)}}_*V_1 = X_{12}^*$, where $X_{12}^*$ denotes the fundamental vector field generated by $X_{12}$. This establishes (\ref{formula_mu_theta}).
\end{proof}

\begin{thm} 
Let $(B, \theta)$ be a $\Gamma$-bundle with connection. For a split $\sigma$ of the Lie algebra, the 1-form $\nabla_{\theta, \sigma}$ on $J_B$ defined by 
\begin{eqnarray}
\nabla_{\theta, \sigma} = 
\zeta^*\nu_\sigma + \rho^*Z_\sigma(\zeta^{-1}, \pi_1^*\theta),
\label{the_connection}
\end{eqnarray}
is a connection on the lifting bundle gerbe $(J_B, B)$, where $\rho : J_B \to B^{[2]}$ is the projection.
\end{thm}

\begin{proof}
This theorem is the consequence of the following two formulas:
\begin{eqnarray}
\pi_{12}^*(\zeta^*\nu) + \pi_{23}^*(\zeta^*\nu) - m^*\pi_{13}^*(\zeta^*\nu)
=
- (\rho_{12} \otimes \rho_{23})^*
Z(\pi_{23}^*\zeta^{-1}, \pi_{12}^*\zeta^*\mu), \\
\pi_{12}^*Z(\zeta^{-1}, \pi_1^*\theta) +
\pi_{23}^*Z(\zeta^{-1}, \pi_1^*\theta) -
\pi_{13}^*Z(\zeta^{-1}, \pi_1^*\theta) =
Z(\pi_{23}^*\zeta^{-1}, \pi_{12}^*\zeta^*\mu),
\end{eqnarray}
where $\rho_{ij} : \pi_{ij}^*J \to B^{[3]}$ is the projection induced by the projection $\rho$ and we omit the subscription of $Z_\sigma$ and $\nu_\sigma$. To prove the first formula we put a tangent vector into both hand sides and compute the values. In the computation we use (\ref{formula_left_nu}) and (\ref{formula_right_nu}). Note also the formula $m_*(v_1, v_2) = {R_{\hgamma_2}}_*v_1 + {L_{\hgamma_1}}_*v_2$, where $m : \hGamma \times \hGamma \to \hGamma$ is the group product and $v_i \in T_{\hgamma_i}\hGamma$ for $i = 1, 2$. The second formula is directly checked by (\ref{formula_Z_Ad}) and (\ref{formula_mu_theta}).
\end{proof}

\begin{prop} \label{prop_conn_diff_sigma}
Let $(B, \theta)$ be a $\Gamma$-bundle with connection. If $\sigma$ and $\sigma'$ are splits of the Lie algebra, then we have
\begin{eqnarray}
\nabla_{\theta, \sigma'} - \nabla_{\theta, \sigma} = 
\rho^* (\pi_1^* - \pi_2^*) \left( (\sigma' - \sigma)(\theta) \right).
\end{eqnarray}
\end{prop}

\begin{proof}
Using $\nu_{\sigma'} - \nu_{\sigma} = -(\sigma' - \sigma)(q^*\mu)$ and (\ref{formula_Z_diff_sigma}), we obtain
$$
\nabla_{\theta, \sigma'} - \nabla_{\theta, \sigma} = 
- (\sigma' - \sigma)(\zeta^*q^*\mu) +
\rho^*(\sigma' - \sigma)(\pi_1^*\theta - \Ad_{\zeta^{-1}}\pi_1^* \theta).
$$
Then the formula follows from $\zeta^*q^*\mu = \rho^*\zeta^*\mu$ and (\ref{formula_mu_theta}).
\end{proof}

In order to construct a curving for the connection on the lifting bundle gerbe we use a reduced splitting with respect to the split $\sigma$ which we choose to define the bundle gerbe connection.

\begin{lem}
{\bf (a)}
For a $\Gamma$-valued function $g$ and a $\Lie\Gamma$-valued 1-form $\alpha$
\begin{eqnarray}
dZ_\sigma(g^{-1}, \alpha) = Z_\sigma(g^{-1}, d\alpha) - \frac{1}{2}
\left\{
\omega_\sigma(g^*\mu, \Ad_{g^{-1}}\alpha) + 
\omega_\sigma(\Ad_{g^{-1}}\alpha, g^*\mu)
\right\}.
\label{formula_dZ1}
\end{eqnarray}

{\bf (b)}
For $\gamma \in \Gamma$ and $X, Y \in \Lie\Gamma$ we have
\begin{eqnarray}
\omega_\sigma(\Ad_\gamma X, \Ad_\gamma Y) = 
\omega_\sigma(X, Y) + Z_\sigma(\gamma, [X, Y]_\Gamma).
\label{formula_ad_omega}
\end{eqnarray}
\end{lem}

\begin{proof}
A direct computation proves (a). For (b), using $Z_\sigma$ and the invariance of the bracket under the adjoint action, we can rewrite the left hand side of (\ref{formula_ad_omega}) as $\Ad_{\gamma}[\sigma(X), \sigma(Y)]_{\hGamma} - \sigma(\Ad_{\gamma}[X, Y]_\Gamma)$. Now the right hand side of (\ref{formula_ad_omega}) follows from $[\sigma(X), \sigma(Y)]_{\hGamma} = \sigma([X, Y]_\Gamma) + \omega_\sigma(X, Y)$.
\end{proof}

\begin{thm}
Let $(B, \theta, L)$ be a $\Gamma$-bundle with connection and splitting. For a split $\sigma$ of the Lie algebra, we define $\kappa_{\theta, \sigma} \in \i\Omega^2(B)$ by 
\begin{eqnarray}
\kappa_{\theta, \sigma}(b; V, W) = \l_\sigma(b, F_\theta(b; V, W)),
\label{the_kappa}
\end{eqnarray}
where $\l_\sigma$ is the reduced splitting with respect to $\sigma$ induced from $L$. Then the curvature of the connection $\nabla_{\theta, \sigma}$ is expressed as
\begin{eqnarray}
F_{\nabla_{\theta, \sigma}} = 
\pi_1^*\left(\frac{1}{2}\omega_\sigma(\theta, \theta) + 
\kappa_{\theta, \sigma} \right) -
\pi_2^*\left(\frac{1}{2}\omega_\sigma(\theta, \theta) + 
\kappa_{\theta, \sigma} \right) .
\end{eqnarray}
Therefore a curving $f_{\theta, \sigma}$ for $\nabla_{\theta, \sigma}$ is given by
\begin{eqnarray}
f_{\theta, \sigma} =
 - \frac{1}{2}\omega_\sigma(\theta, \theta) - \kappa_{\theta, \sigma}.
\label{the_curving}
\end{eqnarray}
\end{thm}

\begin{proof}
This is proved by the following formulas
\begin{eqnarray}
&  &
F_{\nabla_{\theta, \sigma}}  =
\pi_1^*\left(\frac{1}{2}\omega_\sigma(\theta, \theta)\right) -
\pi_2^*\left(\frac{1}{2}\omega_\sigma(\theta, \theta)\right) +
Z_\sigma(\zeta^{-1}, \pi_1^*F_{\theta}) ,\\
&  &
\pi_1^*\kappa_{\theta, \sigma} - \pi_2^*\kappa_{\theta, \sigma} = 
Z_\sigma(\zeta^{-1}, \pi_1^*F_\theta) .
\end{eqnarray}
The first formula follows from (\ref{formula_mu_theta}), (\ref{formula_dZ1}) and (\ref{formula_ad_omega}). To prove the second formula we put tangent vectors into both hand sides and compute the values. In the computation we use (\ref{formula_red_split}) and the property that the value of a curvature is zero for any fundamental vector field.
\end{proof}

\begin{prop}
Let $(B, \theta, L)$ be a $\Gamma$-bundle equipped with connection and splitting. If $\sigma$ and $\sigma'$ are splits of the Lie algebra, then we have 
\begin{eqnarray}
f_{\theta, \sigma'} - f_{\theta, \sigma} = - (\sigma' - \sigma)( d\theta ).
\end{eqnarray}
\end{prop}

\begin{proof}
Easily we have $\omega_{\sigma'}(\theta, \theta) - \omega_{\sigma}(\theta, \theta) = (\sigma' - \sigma)([\theta, \theta]_\Gamma)$. By Corollary \ref{cor_spl_redspl} the reduced splittings $l_\sigma$ and $l_{\sigma'}$ induced from $L$ satisfy (\ref{rel_redspl}). Thus we obtain $\kappa_{\theta, \sigma'} - \kappa_{\theta, \sigma} = (\sigma' - \sigma)(F_\theta)$. These formulas prove the proposition.
\end{proof}

\begin{thm} \label{thm_corres2}
Let $(B, \theta, L)$ be a $\Gamma$-bundle over $X$ with connection and splitting, and $\sigma$ a split of the Lie algebra. In $H^2(X, \F^2)$, the Deligne cohomology class of $(J_B, \nabla_{\theta, \sigma}, f_{\theta, \sigma})$ coincides with the obstruction class of $(B, \theta, L)$.
\end{thm}

\begin{proof}
We show that the \Cech cocycles of the Deligne cohomology classes coincide exactly under an appropriate choice. We fix sections $\{ s_\alpha \}$ of $B$ with respect to a good cover $\U$. First, we take $\h{g}_{\alpha \beta}$ to define $z_{\alpha \beta \gamma}$. Recall the proof of Theorem \ref{thm_corres1}. We can take the section of $(s_\alpha, s_\beta)^*J_B$ as $\h{g}_{\alpha \beta}$. As a consequence $z_{\alpha \beta \gamma} = g_{\alpha \beta \gamma}$ holds. Secondly, to define $u_{\alpha \beta}$ we can take $\h{\theta}_\alpha = \sigma(\theta_\alpha)$. Using $\nu_\sigma$ and $Z_\sigma$ we obtain $u_{\alpha \beta} = - \{ Z_\sigma({g_{\alpha \beta}}^{-1}, \theta_\alpha) + \h{g}_{\alpha \beta}^*\nu_\sigma \} = A_{\alpha \beta}$. Lastly, since $F_{\h{\theta}_\alpha} = \sigma(F_{\theta_\alpha}) + \frac{1}{2}\omega_\sigma(\theta_\alpha, \theta_\alpha)$, we have $K_\alpha = \l_\sigma(s_\alpha, F_{\theta_\alpha}) + \frac{1}{2}\omega_\sigma(\theta_\alpha, \theta_\alpha) = f_\alpha$. 
\end{proof}

This theorem directly establishes the following.

\begin{cor}
The Deligne cohomology class of $(J_B, \nabla_{\theta, \sigma}, f_{\theta, \sigma})$ is independent of the choice of the splits $\sigma$ of the Lie algebra.
\end{cor}

The theorem combined with Theorem \ref{thm_triv_del_bgcc} also gives a corollary which is a generalization of Corollary \ref{cor_top_lbg}.

\begin{cor}
The $(B, \theta, L)$ admits a lifting with a connection compatible with $\theta$ whose scalar curvature vanishes if and only if the induced lifting bundle gerbe $(J_B, \nabla_{\theta, \sigma}, f_{\theta, \sigma})$ is trivial.
\end{cor}

We describe the relation between the curving and the scalar curvature.

\begin{prop} \label{prop_sclcurv_3curv}
Let $(B, \theta, L)$ be a $\Gamma$-bundle with connection and splitting. Suppose that $B$ admits a lifting $(\hB, \hq)$ and $\h{\theta}$ is a connection compatible with $\theta$. If we choose a split $\sigma$ of the Lie algebra, then the curving $f_{\theta, \sigma}$ is expressed as
\begin{eqnarray}
f_{\theta, \sigma} = F_N - \pi^*K_{\h{\theta}},
\label{formula_sclcurv_curving}
\end{eqnarray}
where $F_N$ is the curvature of the connection $N = \h{\theta} - \sigma(\hq^*\theta)$ on the $\T$-bundle $\hq : \hB \to B$. Moreover, the 3-curvature of $(J_B, \nabla_{\theta, \sigma}, f_{\theta, \sigma})$ is expressed as
\begin{eqnarray}
\Xi = - d K_{\h{\theta}} .
\label{formula_sclcurv_3curv}
\end{eqnarray}
\end{prop}

\begin{proof}
If $\h{\theta}$ is a compatible with $\theta$, then the 1-form $N$ takes values in $\i\R$. The basic properties of connection show that $N$ is indeed a connection. It is straight to see $F_{\h{\theta}} = F_{N} + \hq^*(\sigma(F_\theta) + \frac{1}{2}\omega_\sigma(\theta, \theta))$. Note that we can regard $F_{\h{\theta}}$ as a 2-form on $B$. Using this formula and the reduced splitting $\l_\sigma$, we have
\begin{eqnarray*}
\pi^*K(b; V, W) 
& = & 
\l_\sigma(b, F_\theta(b; V, W)) + 
\left( F_N + \frac{1}{2}\omega_\sigma(\theta, \theta) \right) (b; V, W) \\
& = & 
\left( F_N - f_{\theta, \sigma} \right) (b; V, W)
\end{eqnarray*}
for $V, W \in T_bB$. Differentiating (\ref{formula_sclcurv_curving}), we have (\ref{formula_sclcurv_3curv}).
\end{proof}

Let $B$ be a $\Gamma$-bundle over $X$, and $\varphi : X' \to X$ a map. By the construction of the adjoint bundle, we have a natural isomorphism $\varphi^*(B \times_{Ad} \Lie\hGamma) \cong (\varphi^*B) \times_{Ad} \Lie\hGamma$ of vector bundles. Hence a splitting $L$ of $B$ induces a splitting of $\varphi^*B$. We denote the splitting by $\varphi^*L$ and call it the \textit{pull-back} of $L$.

\begin{thm}
Let $(B, \theta, L)$ be a $\Gamma$-bundle over $X$ with connection and splitting, $\sigma$ a split of the Lie algebra, and $\varphi : X' \to X$ a map. We denote the lifting bundle gerbes with connection and curving associated with $(B, \theta, L)$ and with $\varphi^*(B, \theta, L)$ by $(J_B, \nabla_{\theta, \sigma}, f_{\theta, \sigma})$ and by $(J_{\varphi^*B}, \nabla_{\varphi^*\theta, \sigma}, f_{\varphi^*\theta, \sigma})$ respectively. Then there is a natural isomorphism $\varphi^*(J_B, \nabla_{\theta, \sigma}, f_{\theta, \sigma}) \cong (J_{\varphi^*B}, \nabla_{\varphi^*\theta, \sigma}, f_{\varphi^*\theta, \sigma})$.
\end{thm}

\begin{proof}
The lifting bundle gerbe of the pull-back is by definition $J_{\varphi^*B} = {\zeta'}^*\hGamma$, where $\zeta' : (\varphi^*B)^{[2]} \to \Gamma$ is defined in the same way as $\zeta$. If $\til{\varphi} : \varphi^*B \to B$ is the bundle map covering $\varphi$, then we have $\zeta' = \zeta \circ \til{\varphi}^{[2]}$. Using this relation, we obtain a natural isomorphism $\varphi^*(J_B, B) \cong (J_{\varphi^*B}, \varphi^*B)$ of bundle gerbes. Because the pull-back of the connection $\theta$ is $\til{\varphi}^*\theta$, we can see that the isomorphism maps the connection on $\varphi^*(J_B, B)$ to that on $(J_{\varphi^*B}, \varphi^*B)$ by (\ref{the_connection}). If the reduced splitting $\l_\sigma : B \times \Lie\Gamma \to \i\R$ is induced by the splitting $L$, then the reduced splitting induced by the pull-back $\varphi^*L$ is $\l_\sigma \circ (\til{\varphi} \times \id) : \varphi^*B \times \Lie\Gamma \to \i\R$. So we obtain $\til{\varphi}^*f_{\theta, \sigma} = f_{\til{\varphi}^*\theta, \sigma}$ by (\ref{the_kappa}) and (\ref{the_curving}).
\end{proof}


\section{Examples}
\label{sec_example}

We give an important example of lifting problems. We set $G = SU(2)$ and denote by $LG = \M{S^1}{G}$ the space of loops in $G$.  By the pointwise product $LG$ is an infinite dimensional Lie group known as the \textit{loop group} \cite{P-S}. The Lie algebra of the loop group is $L\g$, the space of loops in the Lie algebra $\g$ of $G$. It is well-known that there exists a central extension $\h{LG}$ of $LG$ called the \textit{Kac-Moody group} \cite{P-S} for each $k \in \Z$. We can take a split of $\h{L\g}$ such that the Lie algebra cocycle is given by
\begin{eqnarray}
\omega(X, Y) = \frac{k\i}{2\pi} \int_{S^1} \Tr (X dY).
\label{formula_cocycle_KM}
\end{eqnarray}
We fix the split and omit subscriptions.

\begin{prop}
The map $Z : LG \times L\g \to \i\R$ introduced in Definition \ref{dfn_group_cocycle} is given by
\begin{eqnarray}
Z(\gamma, X) = -\frac{k\i}{2\pi}\int_{S^1} 
\Tr \left( \gamma^{-1}d\gamma X \right).
\end{eqnarray}
\end{prop}

\begin{proof}
Using (\ref{formula_ad_omega}) and (\ref{formula_cocycle_KM}), we have
\begin{eqnarray}
Z(\gamma, [X, Y]) = -\frac{k\i}{2\pi}\int_{S^1} 
\Tr \left( \gamma^{-1}d\gamma \ [X, Y] \right).
\end{eqnarray}
Because $SU(2)$ is simple, we have $\g = [\g, \g]$. This implies that $L\g \subset [L\g, L\g]$, and the proof is completed.
\end{proof}

We can also obtain this formula by the computation based on an explicit construction of the central extension \cite{Mi}.

If we consider a principal $G$-bundle $\pi : P \to M$, then we have a natural $LG$-bundle $\pi_L : LP \to LM$ called the \textit{loop bundle}. In this case, a lifting is called a \textit{string structure} and the obstruction class is called the \textit{string class} \cite{C-P, K}. 

It is easy to see that a connection $A$ on $P$ induces a connection $\bar{A}$ on $LP$ by $\bar{A}(p; V, W) = A(p(t); V(t), W(t))$. Moreover, the connection $A$ gives a reduced splitting.

\begin{prop}
A connection $A$ on $P$ induces a reduced splitting of $LP$ by
\begin{eqnarray}
\l(p,X) = - \frac{k\i}{2\pi} \int_{S^1} \Tr \left( p^*A \ X \right).
\label{formula_explict_redspl}
\end{eqnarray}
\end{prop}

\begin{proof}
It is clear that $\l$ is linear if $p \in LP$ is fixed. So we show that $\l$ satisfies (\ref{formula_red_split}). We have the following relation between $\g$-valued 1-forms on $S^1$.
\begin{eqnarray}
(p\gamma)^*A = \Ad_{\gamma^{-1}} p^*A + \gamma^{-1} d\gamma.
\label{formula_forms_on_circle}
\end{eqnarray}
Since $\Tr$ is invariant under the adjoin action of $G$, the proposition is proved.
\end{proof}

By the result of previous section we obtain the lifting bundle gerbe with connection and curving $(J_{LP}, \nabla_{\bar{A}}, f_{\bar{A}})$. The curving is explicitly given as 
\begin{eqnarray}
&   &
f_{\bar{A}}(p; V, W) \label{formula_explicit_curv} \\
&   & = 
- \frac{k\i}{2\pi} \int_{S^1}
\left\{
\frac{1}{2} \Tr \left( A(p;V) d(A(p; W)) \right)
-
\Tr \left( p^*A \ F_A(p; V, W) \right)
\right\}. \nonumber
\end{eqnarray}

This 2-form on $LP$ appears in the work of Coquereaux and Pilch \cite{C-P}. To explain their work, we introduce the \textit{transgression maps} or \textit{averaged evaluations} for differential forms. Let $F$ be a compact oriented $m$-dimensional manifold, and $\alpha$ an $n$-form on a manifold $X$. We denote by $\M{F}{X}$ the space of smooth maps from $F$ to $X$. We pull $\alpha$ back by the evaluation map $\ev : \M{F}{X} \times F \to X$, and apply the fiber integration to $\M{F}{X} \times F \to \M{F}{X}$. Then we have an $(n-m)$-form $\tau_F\alpha = \int_F\ev^*\alpha$ on $\M{F}{X}$ and the transgression map $\tau_F : \Omega^n(X) \to \Omega^{n-m}(\M{F}{X})$. The transgression and the exterior derivative satisfy Stokes' formula $(-1)^m d\tau_F\alpha = \tau_Fd\alpha - \tau_{\partial F}\alpha$. In the case that $F = S^1$ we write the transgression map as $\tau_L$.

There is a relation between 2-forms on $LG$
\begin{eqnarray}
-\frac{1}{2} \omega(\mu, \mu) = 
2\pi\i \left( \tau_L\sigma + d\beta \right),
\label{formula_volume_G}
\end{eqnarray}
where $\sigma \in \Omega^3(G)$ and $\beta \in \Omega^1(LG)$ are given by
\begin{eqnarray}
\sigma & = & \frac{k}{24\pi^2} \Tr (g^{-1}dg )^3, \\
\beta(\gamma; {L_\gamma}_*X) & = &
\frac{k}{8\pi^2} \int_{S^1} \Tr(\gamma^{-1}d\gamma X).
\end{eqnarray}

In \cite{C-P}, Coquereaux and Pilch described the following relation as a fiberwise extension of (\ref{formula_volume_G}):
\begin{eqnarray}
f_{\bar{A}} = 2\pi\i \left( \tau_L CS + d \Upsilon \right) ,
\label{formula_curv_CS}
\end{eqnarray}
where $CS \in \Omega^3(P)$ is the Chern-Simons form \cite{C-S}
\begin{eqnarray}
CS = \frac{k}{8\pi^2} 
\Tr \left( A \wedge dA + \frac{2}{3}A \wedge A \wedge A \right)
\end{eqnarray}
and $\Upsilon \in \Omega^1(LP)$ is defined by
\begin{eqnarray}
\Upsilon(p; V) = 
\frac{k}{8\pi^2} \int_{S^1} \Tr \left( p^*A \ A(p; V) \right).
\end{eqnarray}

Using the formula (\ref{formula_curv_CS}), we can see that the 3-curvature is $\Xi = -2\pi\i \tau_L\xi$, where $\xi = \frac{1}{8\pi^2} \Tr(F_A \wedge F_A)$ is the characteristic 4-form on $M$ induced by $A$. From this we find the fact known in \cite{C-P} that a de Rham representative of the string class of the loop bundle is given by the transgression of the 4-form.

\medskip

Now we take an oriented compact 2-manifold $\Sigma$ with boundary $\partial \Sigma = S^1$. We denote the restriction map by $r : \M{\Sigma}{M} \to LM$. On the pull-back $r^*LP$ we can construct a string structure as follows. The map space $\G = \M{\Sigma}{G}$ is a Lie group and $\pi_\Sigma : \M{\Sigma}{P} \to \M{\Sigma}{M}$ is a principal $\G$-bundle. We define a map $c : \M{\Sigma}{P} \times \G \to \h{LG}$ by
\begin{eqnarray}
c(f, g) = \exp \frac{k\i}{4\pi} \int_\Sigma 
\Tr (f^*A \wedge dg g^{-1}) \cdot \left( e^{W_\Sigma(g)} \right)^{-1},
\end{eqnarray}
where $e^{W_\Sigma(g)} \in \h{LG}$ is the \textit{Wess-Zumino term} \cite{F}. If $q : \h{LG} \to LG$ is the projection, then we have $q(c(f, g)) = \partial g^{-1}$.

\begin{prop}
We define a quotient space $\hB_\Sigma = \M{\Sigma}{P} \times \h{LG} / \sim$ by the equivalence relation $(f, \hgamma) \sim (fg, c(f, g)\hgamma)$ for $g \in \G$. If we define a right action of $\h{LG}$ by $[f, \hgamma] \cdot \h{\eta} = [f, \hgamma \cdot \h{\eta}]$ and a projection $\hB_\Sigma \to \M{\Sigma}{M}$ by $[f, \hgamma] \mapsto \pi_\Sigma(f)$, then $\hB$ is a principal $\h{LG}$-bundle. Moreover, if we define $\hq_\Sigma : \hB_\Sigma \to r^*LP$ by $\hq_\Sigma([f, \hgamma]) = (\pi_\Sigma(f), \partial f \cdot q(\hgamma))$, then $(\hB_\Sigma, \hq_\Sigma)$ is a string structure on $r^*LP$. 
\end{prop}

\begin{proof}
We can check $c(f, gh) = c(fg, h) c(f, g)$ by the Polyakov-Wiegmann formula of the Wess-Zumino term \cite{F}. Thus the equivalence relation is well-defined. Since the action of $\G$ on $\M{\Sigma}{P}$ is free, the fiber of $\hB_\Sigma \to \M{\Sigma}{M}$ is the orbit of the right action. Therefore we obtain a principal $\h{LG}$-bundle. It is easy to see that $\hq_\Sigma$ is an equivariant map, and we obtain a string structure.
\end{proof}

The existence of string structure implies that the Dixmier-Douady class $d(J_{r^*LP}) = r^*d(J_{LP})$ is trivial and the 3-curvature $r^*\Xi$ is exact. If we have a connection on $\hB_\Sigma$ compatible with $r^*\bar{A}$, then we can compute a 2-form cobounding $r^*\Xi$ by means of Proposition \ref{prop_sclcurv_3curv}.

\begin{prop}
We define a 1-form on $\M{\Sigma}{P}$ with values in $\i\R$ by 
\begin{eqnarray}
\til{N} = -2\pi\i(\tau_\Sigma CS - r^*\Upsilon), 
\end{eqnarray}
where $r : \M{\Sigma}{P} \to LP$ is the restriction map. Then the 1-form $\hgamma^{-1} (r^*\bar{A} + \til{N}) \hgamma + \h{\mu}$ on $\M{\Sigma}{P} \times \h{LG}$ descends to a connection $\h{A}$ on $\hB_\Sigma$ compatible with $r^*\bar{A}$ whose scalar curvature is
\begin{eqnarray}
K_{\h{A}} = -2\pi\i \tau_\Sigma \xi .
\end{eqnarray}
\end{prop}

\begin{proof}
Note that there is the following commutative diagram.
$$
\begin{CD}
\M{\Sigma}{P} \times \h{LG} @>>> 
\M{\Sigma}{P} \times LG @>>> 
\M{\Sigma}{P} \\
@VVV @VVV @VVV \\
\hB_\Sigma @>>> r^*LP @>>> \M{\Sigma}{M}
\end{CD}
$$
Here the lows are string structures and the columns are principal $\G$-bundles.
Thus, we can investigate the data of the lower string structure by that of the upper string structure. It is obvious that the connection $\hgamma^{-1} (r^*\bar{A} + \til{N}) \hgamma + \h{\mu}$ on $\M{\Sigma}{P} \times \h{LG}$ is compatible with the connection $\gamma^{-1}(r^*\bar{A})\gamma + \mu$ on $\M{\Sigma}{P} \times LG$. By computation we can show that these connections descend. In particular, the latter connection descends to $r^*\bar{A}$. Hence we obtain the compatible connection $\h{A}$ on $\hB_\Sigma$. For the computation of the scalar curvature, we also use the commutative diagram. By Proposition \ref{prop_sclcurv_3curv} and the formula
\begin{eqnarray}
\pi_\Sigma^*(-2\pi\i \tau_\Sigma \xi) =
d\til{N} - r^*f_{\bar{A}}
\end{eqnarray}
derived from Stokes' formula, we obtain the result.
\end{proof}

We can check that $r^*\Xi = -dK_{\h{A}}$ by using Stokes' formula.

\medskip

Recently, Murray and Stevenson \cite{Mu-S2} also described connections and curvings on the lifting bundle gerbes associated with loop group bundles. In the case of $(LP, \bar{A})$ we can verify that $\nabla_{\bar{A}}$ coincides with their description. While the curving in this paper is given by using a reduced splitting, the curving in \cite{Mu-S2} is given by using a \textit{twisted Higgs field}.

\begin{dfn}[\cite{Mu-S2}]
A \textit{twisted Higgs field} $\Phi$ on a principal $LG$-bundle $B$ is a function $\Phi : B \to \M{[0,2\pi]}{\g}$ such that
\begin{eqnarray}
\Phi_{b\gamma} = \Ad_{\gamma^{-1}} \Phi_{b} 
+ \gamma^{-1} \frac{\partial \gamma}{\partial \theta},
\label{formula_thf}
\end{eqnarray}
where $\Phi_b : [0, 2\pi] \rightarrow \g$ is the image of $b \in B$.
\end{dfn}

On $LP$ we have a twisted Higgs field induced from a connection on $P$.

\begin{prop}
Let $B$ be a loop bundle $LP$. For $p \in LP$ we define $\Phi_p : [0, 2\pi] \to \g$ by $p^*A = \Phi_p(\theta) d\theta$, where $S^1$ is identified with $[0, 2\pi]$. Then $\Phi$ is a twisted Higgs field on $LP$.
\end{prop}

\begin{proof}
If we use the coordinate $\theta \in S^1$, then the relation (\ref{formula_forms_on_circle}) is equivalent to the relation (\ref{formula_thf}).
\end{proof}

When we express (\ref{formula_explict_redspl}) and (\ref{formula_explicit_curv}) using the coordinate $\theta \in S^1$, we see that the reduced splitting and the curving include the twisted Higgs field. If we use this twisted Higgs field in the result of \cite{Mu-S2}, then we have the same curving as (\ref{formula_explicit_curv}).

\bigskip

\textit{Acknowledgments}. I would like to thank M.K. Murray and D. Stevenson for sending their preprint and for useful comments. I would like to thank T. Kohno for helpful suggestions and for constant encouragements. 



Graduate school of Mathematical Sciences, 
University of Tokyo, Komaba 3-8-1, Meguro-Ku, Tokyo, 153-8914 Japan.

e-mail: kgomi@ms.u-tokyo.ac.jp


\begin{thebibliography}{999}

\bibitem{B-T}R. Bott and L. Tu,
\textit{Differential forms in algebraic topology},
Graduate Texts in Mathematics, 82, Springer-Verlag, New York-Berlin, 1982.

\bibitem{Bro}K. S. Brown,
\textit{Cohomology of groups},
Graduate Texts in Mathematics, 87, Springer-Verlag, New York-Berlin, 1982. 

\bibitem{B}J-L. Brylinski,
\textit{Loop spaces, characteristic classes and geometric quantization},
Progress in Mathematics, 107, Birkh$\ddot{\textrm{a}}$user Boston, Inc., Boston, MA, 1993.

\bibitem{B-M}J-L. Brylinski and D. A. McLaughlin,
`The geometry of degree-four characteristic classes and of line bundles on loop spaces. I',
Duke. Math. J. 75 (1994), no.3, 603-638.

\bibitem{C-P}R. Coquereaux and K. Pilch,
`String structures on loop bundles',
Comm. Math. Phys. 120 (1989), no.3, 353-378.

\bibitem{C-S}S. S. Chern and J. Simons,
`Characteristic forms and geometric invariants',
Ann. Math. 99 (1974), 48-69.

\bibitem{F}D. S. Freed,
`Classical Chern-Simons Theory. I',
Adv. Math. 113 (1995), no.2, 237-303.

\bibitem{Gi}J. Giraud,
\textit{Cohomologie  non-ab\'elienne},
Die Grundlehren der mathematischen Wissenschaften, Band 179. 
Springer-Verlag, Berlin-New York, 1971.

\bibitem{K}T. P. Killingback,
`World-sheet anomalies and loop geometry',
Nucl. Phys. B 288 (1987), no.3-4, 578-588.

\bibitem{Mi}J. Mickelsson,
`Kac-Moody groups and the Dirac determinant line bundle',
\textit{Topological and geometrical methods in field theory} (Espoo, 1986), Edited by J. Hietarinta and J. Westerholm. World Sci. Publishing, Teaneck, NJ, 1986, pp.117-131.

\bibitem{Mu}M. K. Murray,
`Bundle gerbes',
J. London Math. Soc. (2) 54 (1996), no.2, 403-416.

\bibitem{Mu-S1}M. K. Murray and D. Stevenson,
`Bundle gerbes: stable isomorphism and local theory',
J. London Math. Soc. (2) 62 (2000), no.3, 925-937.

\bibitem{Mu-S2}M. K. Murray and D. Stevenson,
`Higgs fields, bundle gerbes and string structures',
math.DG/0106179.

\bibitem{P-S}A. Pressley and G. Segal,
\textit{Loop groups}, Oxford Mathematical Monographs, Oxford University Press, New York, 1986.
\end{thebibliography}
\end{document}